\documentclass[11pt]{article}
\usepackage[latin1]{inputenc}
\usepackage{amssymb}
\usepackage{latexsym}
\usepackage{color}
\usepackage{amsmath}
\usepackage{amsfonts}
\usepackage{amsthm}
\usepackage{varioref}
\usepackage{enumerate,pdfsync}
\usepackage{cancel}

\def\be{\begin{eqnarray}}
\def\ee{\end{eqnarray}}
\def\beq{\begin{equation}}
\def\eeq{\end{equation}}
\def\b*{\begin{eqnarray*}}
\def\e*{\end{eqnarray*}}
\def\bi{\begin{itemize}}
\def\ei{\end{itemize}}

\def \1{{\bf 1}}

\def\eps{\varepsilon}

\def\={\;=\;}

\def\x{\times}

\DeclareMathOperator*{\esssup}{ess\,sup}

\def \proof{{\noindent \bf Proof. }}
\def \ep{\hbox{ }\hfill$\Box$}
 \def\reff#1{{\rm(\ref{#1})}}

 \def\vs#1{\vspace{#1mm}}

\def \E{\mathbb{E}}
\def \F{\mathbb{F}}

\def \P{\mathbb{P}}
\def \Q{\mathbb{Q}}
\def \R{\mathbb{R}}

\def\Ac{{\cal A}}

\def\Dc{{\cal D}}
\def\Ec{{\cal E}}
\def\Fc{{\cal F}}
\def\Gc{{\cal G}}

\def\Kc{{\cal K}}

\def\Tc{{\cal T}}
\def\Uc{{\cal U}}

\def\Yc{{\cal Y}}
\def\Zc{{\cal Z}}

\def\Lb{{\mathbf L}}
\def \Hb{{\mathbf H}}

\def\Yb{\bar Y}
\def\Zb{\bar Z}

\newtheorem{Theorem}{Theorem}[part]

\newtheorem{Proposition}{Proposition}[part]

\newtheorem{Lemma}{Lemma}[part]
\newtheorem{Corollary}{Corollary}[part]
\newtheorem{Remark}{Remark}[part]

\makeatletter \@addtoreset{equation}{section}

\@addtoreset{Definition}{section}

\@addtoreset{Theorem}{section}

\@addtoreset{Proposition}{section}

\@addtoreset{Assumption}{section}

\@addtoreset{Corollary}{section}

\@addtoreset{Lemma}{section}

\@addtoreset{Remark}{section}

\@addtoreset{Example}{section}

\def\abs#1{\left|#1\right|}

\theoremstyle{plain}
\newtheorem{theorem}{Theorem}[section]
\newtheorem{proposition}[theorem]{Proposition}

\theoremstyle{definition}
\newtheorem{remark}[theorem]{Remark}

\def\Sb{{\bf S}}
\def\Kb{{\bf K}}

\def\BSDE{{\rm BSDE}}

\def\Yb {{\tilde{\cal Y}}}
\def\Zb {{\tilde{\cal Z}}}
\def\Kb {{\tilde{\cal K}}}

\def\Ybb {\hat \Yc}
\def\Zbb {\hat \Zc}
\def\Kbb {\hat \Kc}

\def\Ecf#1#2#3{\tilde \Ec_{#1}^{#2}[#3]}

 \def\Db{\mathbf{\mathbf D}}
 \def\abs#1{\left|#1\right|}

\def\bal{\begin{aligned}}
\def\eal{\end{aligned}}

\renewcommand{\thefootnote}{\*}

\def\BSDE{{\rm BSDE}}
\def\thetar{{\mathrm \theta}}
\def\[{[\![}
\def\]{]\!]}

\title{ Regularity of BSDEs with a convex constraint on the gains-process }
\author{Bruno Bouchard\footnote{Universit\'e Paris-Dauphine CEREMADE UMR CNRS 7534 $\&$ CREST, email \texttt{bouchard@ceremade.dauphine.fr}} \and Romuald Elie\footnote{Universit\'e Paris-Est, LAMA UMR CNRS 8050, email \texttt{romuald.elie@univ-mlv.fr}} \and Ludovic Moreau\footnote{ETH Z\"urich, Departement f\"ur Mathematik, email \texttt{ludovic.moreau@math.ethz.ch}. Partially supported by the Swiss National Science Foundation under the grant SNF $200021\_152555$ and by the ETH Foundation.} \and \\
}

\begin{document}
\maketitle
\begin{abstract}  We\footnote{Research supported by ANR Liquirisk and Investissements d'Avenir (ANR-11-IDEX-0003/Labex Ecodec/ANR-11-LABX-0047).} {consider} the minimal super-solution  of a backward stochastic differential equation with constraint on the gains-process. The terminal condition is given by a function of the terminal value of a forward stochastic differential equation. Under boundedness assumptions on the coefficients, we show that the first component of the solution is Lipschitz in space and $\frac12$-H\"older in time with respect to the initial data of the forward process. Its path is continuous before the time horizon at which its left-limit is given by a face-lifted version of its natural boundary condition. This first component is actually equal to its own face-lift. We only use probabilistic arguments. In particular, our results can be extended to certain non-Markovian settings.

\end{abstract}
\renewcommand{\thefootnote}{\arabic{footnote}} 
\setcounter{footnote}{0}
\vspace{5mm}


\vspace{5mm}

\noindent {\textbf{MSC Classification (2010)}: 60H10, 60H30, 49L20.

\section{Introduction}

The aim of this paper is to establish new stability results for the minimal super-solution $(\hat \Yc^{\zeta},\hat \Zc^{\zeta})$ of a  backward differential equation of the form
\b*
U_{t}&=& g(X^{\zeta}_{T})+\int_{t}^{T} f(X^{\zeta}_{s},U_{s},V_{s})ds -\int_{t}^{T} V_{s} dW_{s},
\;t\le T,
\e*
satisfying the constraint
\b*
\hat \Zc^{\zeta}\sigma(X^{\zeta})^{-1} \in K\;\;dt\otimes d\P{\rm -a.e.}
\e*
In the above, $W$ is a $d$-dimensional Brownian motion and   $X^{\zeta}$ solves a forward stochastic differential equation with volatility parameter $\sigma$, indexed by the initial conditions $\zeta=(t,x)\in [0,T]\x \R^{d}$: $X^{\zeta}_{t}=x$.

Estimates on the regularity can be of important use in many applications, in particular in the design of probabilistic numerical schemes which,  to the best of our knowledge, are missing for such constrained backward   differential equations. 

When $K=\R^{d}$, i.e. there is no constraint, and the coefficients are Lipschitz continuous,  it is well-known that   $\hat \Yc^{\zeta}$ has continuous path and that the (deterministic) map $(t,x)\mapsto \hat \Yc^{(t,x)}_{t}$ is $1/2$-H\"older in time and Lipschitz in space:
\be\label{eq: intro regu}
|\hat \Yc^{(t,x)}_{t}-\hat \Yc^{(t',x')}_{t'}|\le C\; (|t-t'|^{\frac12}+|x-x'|).
\ee
 See e.g.~\cite{Par98}. This basically follows from standard estimates using It\^{o}'s and Gronwall's Lemma.

In the general case, such a minimal super-solution solves an equation of the form
\b*
\hat \Yc^{\zeta}_{t}&=& g(X^{\zeta}_{T})+\int_{t}^{T} f(X^{\zeta}_{s},\hat \Yc^{\zeta}_{s},\hat \Zc^{\zeta}_{s})ds -\int_{t}^{T} \hat \Zc^{\zeta}_{s} dW_{s}+ \hat \Kc^{\zeta}_{T}-\hat \Kc^{\zeta}_{t},
\;t\le T,
\e*
in which $\hat \Kc^{\zeta}$ is an adapted non-decreasing process, see \cite{CKS,peng1999monotonic}. Because little is known on the regularity of this process, the technics used in the unconstrained case can not  be reproduced.

Nevertheless, it is well-known that such a minimal super-solution can be approximated by a sequence of penalized unconstrained stochastic backward differential equations, see \cite{CKS,peng1999monotonic}. It is therefore tempting to use the  estimates associated to each element of the approximating sequence and to pass to the limit. Unfortunately, the Lipschitz continuity coefficients of the approximating sequence blow-up.

Another way to proceed consists in using the dual formulation of \cite{cvitanic1993hedging,CKS}. In their representation, the component $\hat \Yc^{\zeta}$ is identified to the value of the optimal control problem of a family of backward stochastic differential equations written under a suitable set of equivalent probability measures, see Section \ref{sec: def weak form}. The main difficulty is that it is singular: each of the controls is bounded, but the bound is not uniform.

In this paper, we essentially make use of this dual formulation, but we use a strong version:  the controls are directly incorporated in the dynamics rather than through changes of measures. See Section \ref{sec: estimates via the strong dual formulation}.  Space stability is essentially obvious for this strong version, while it is not in its original 'weak' form. Still, the singularity of the optimal control problem makes estimates on the stability in time quite delicate a-priori.

The key idea of this paper is to use the fact that the solution is automatically 'face-lifted'  in the sense of Proposition \ref{prop: Yb self face lifted}  below. This 'face-lifting' phenomenon is well-known as far as the terminal condition is concerned, this goes back to \cite{Broadie} in the specific setting of mathematical finance, see also \cite{bouchard2010obstacle,cvitanic1999super} and the references therein. We use probabilistic arguments to show that it holds also in the parabolic interior $[0,T)\x \R^{d}$ of the domain, which can be guessed in the setting of \cite{bouchard2010obstacle,cvitanic1999super}  from their pde characterization. This 'face-lifting'  phenomenon allows to absorb the singular control, and to extend \reff{eq: intro regu} to the constrained case. See Theorem \ref{thm: regularity constrained bsde}.
\\

The paper is organized as follows. The setting and the main results are stated in Section \ref{sec: main result}. Comments on our assumptions and possible extensions are discussed at the end in Section \ref{sec: possible extensions}. Section \ref{sec: estimates via the strong dual formulation} is dedicated to the strong version of the dual formulation for which our estimates are established. This part contains the main ideas of this paper. In Section \ref{sec: weak versus strong}, we show that the strong dual formulation coincides with its weak version,  and that the latter actually provides (as well-known  {when the driver is convex}) the first component $\hat \Yc^{\zeta}$ of the constrained backward stochastic differential equation.
\\

{\bf Notations: }
All over this paper, we let $\Omega=C([0,T],\R^{d})$, $d\ge 1$, $T>0$, be the canonical space of continuous $d$-dimensional functions $\omega$ on $[0,T]$ such that $\omega_{0}=0$. It is  endowed with the Wiener measure $\P$. We let $W$ be the coordinate process, $W_{t}(\omega)=\omega_{t}$, and we denote by $\F=(\Fc_{t})_{t\le T}$ the augmentation of its raw filtration under $\P$. Random variables are defined  on $(\Omega,\Fc_{T},\P)$.
For the expectation under $\P$, we simply use the symbol $\E$, while we write $\E^{\Q}$ if it is taken under a different measure $\Q$.
Given a probability measure $\Q$, $\Gc\subset \Fc_{T}$, $p\ge 1$ and $A\subset \R^{n}$, we denote by $\Lb_{p}(A,\Q,\Gc)$ the set of $\Gc$-measurable $A$-valued random variables whose $p$-moment under $\Q$ is finite. We let $\Sb_{2}(\Q)$ (resp. $\Hb_{2}(\Q)$)be the set of $\R^{n}$-valued progressively measurable processes $V$ such that $\E^{\Q}[\sup_{t\le T}|V_{t}|^{2}]<\infty$ (resp. $\E^{\Q}[\int_{t}^{T}|V_{t}|^{2}dt]<\infty$),  in which $|v|$ denotes the Euclydian norm of $v\in \R^{n}$ and $n$ is given by the context.  The set of stopping times with values in $[0,T]$ is $\Tc$, while $\Tc_{\tau}$ is the set of stopping times a.s.~greater than $\tau\in \Tc$. Finally,  $\Db_{2}(\Q)$ denotes the set of couples $(\tau,\xi)\in \Tc\x\Lb_{2}(\Q,\R^{d}, \Fc_{T})$ such that $\xi$ is $\Fc_{\tau}$-measurable.  For $\zeta\in \Db_{2}$, we write $(\tau_{\zeta},\xi_{\zeta})=\zeta$. In all these definitions, we omit the arguments that can be clearly identified by the context. When nothing else is specified, inequalities between random variables or convergence of sequences of random variables hold in the $\P$-a.s.~sense.

\section{Main regularity and stability results}\label{sec: main result}

 As a first step of analysis, we   concentrate on a Markovian setting with rather stringent boundedness assumptions. Possible extensions will be discussed in Section \ref{sec: possible extensions} at the end of this paper. They include another type of constraint, certain non-Markovian settings and optimal control problems.
 \vs2

The forward component is  the unique strong solution $X^{\zeta}$ on $[0,T]$ of the stochastic differential equation
\beq\label{eq: SDE X}
    X^{\zeta}_{t\vee \tau_{\zeta}}=\xi_{\zeta}+\int_{\tau_{\zeta}}^{t\vee \tau_{\zeta}} b_{s}(X^{\zeta}_{s})ds + \int_{\tau_{\zeta}}^{t\vee \tau_{\zeta}} \sigma_{s}(X^{\zeta}_{s})dW_{s},
\eeq
in which the initial data    $\zeta\in \Db_{2}$, and  the parameters $(b,\sigma): [0,T]\x \R^{d}\mapsto \R^{d}\x \R^{d\times d}$ are measurable maps. They are  assumed to be bounded, and Lipschitz in their space variable, uniformly in their time argument.  We also assume that $\sigma$ is invertible with bounded inverse. Namely, there exists $L>0$ such that
   \be\label{eq: ass b sigma}
        |(b_{t},\sigma_{t})(x)-(b_{t},\sigma_{t})(x')|\le L|x-x'|&\mbox{ and } &   (|b_{t}|+|\sigma_{t}|+|\sigma^{-1}_{t}|)(x)\le L,
   \ee
  for all  $(t,x,x')\in [0,T]\x \R^{d}\x \R^{d}$.

The backward equation is defined by  two measurable maps $f: [0,T]\x\R^{d}\x \R \x \R^{d}\mapsto \R$ and $g:\R^{d}\mapsto \R$ such that,    for all $(t,x,\theta), (t,x',\theta')\in  [0,T]\x\R^{d}\x(\R\x \R^{d})$,
 \be\label{eq: ass f}
       &\abs{f_{t}(x,\theta)-f_{t}(x',\theta')}\le L \left(|x-x'|+|\theta-\theta'|\right)\;,\;
        |f_{t}(x,\theta)|\le L(1+|\theta|) \;\;,&
        \\
             &    |g(x)|\le L,\; \mbox{ and } g\mbox{ is lower-semicontinuous.}&\label{eq: ass g}
    \ee
A supersolution of  \BSDE$(f,g,\zeta)$ is a  process $(U,V) \in \Sb_{2}\x \Hb_{2}$ satisfying
\be\label{eq: def constrained bsde 1}
\begin{array}{rcl}
U_{t\vee \tau_{\zeta}}&\ge& U_{t'\vee \tau_{\zeta}}+\int_{t\vee \tau_{\zeta}}^{t'\vee \tau_{\zeta}} f_{s}(X^{\zeta}_{s},U_{s},V_{s})ds -\int_{t\vee \tau_{\zeta}}^{t'\vee \tau_{\zeta}} V_{s} dW_{s},
\;t\le t'\le T.\\
U_{T}&=&g(X^{\zeta}_{T}).
\end{array}
\ee
The constraint on the $V$ coordinate is associated to   a  family  $(K_{t})_{t\le T}$ of closed convex  sets of $\R^{d}$:
\beq\label{eq: def constrained bsde 2}
V\sigma(X^{\zeta})^{-1} \in K\;\;dt\otimes d\P{\rm -a.e} \; \mbox{ on } [\![\tau_{\zeta},T]\!].
\eeq
When it is satisfied, we say that  $(U,V)$ is a super-solution of  BSDE$_{K}(f,g,\zeta)$. This super-solution is said to be minimal if $U'_{s\vee \tau_{\zeta}}\ge U_{s\vee \tau_{\zeta}}$ for all $s\le T$  and any other super-solution $(U',V')\in \Sb_{2}\x \Hb_{2}$ of \BSDE$_{K}(f,g,\zeta)$.

 We require that
  \be
    &0\in  K_{t}\mbox{ for all $t\le T$}&\label{eq: ass K int 0}
         \\
    & \cup_{t\le T} K_{t} \mbox{ is bounded,}& \label{eq: ass K borne}
  \ee
  and that, for all  $u\in \R^{d}$,
  \be
     &t\in [0,T]\mapsto \delta_{t}(u):=\sup \{\;k^{\top} u , k\in K_{t} \}\;\mbox{ is   left-continuous at $T$,} &   \label{eq : ass delta continuous}
      \\
      &  \mbox{ and non-increasing.}&    \label{eq : ass delta non-increasing}
     \ee
\begin{Remark} The conditions \reff{eq: ass K borne}-\reff{eq : ass delta non-increasing} are equivalent to : the family $(K_{t})_{t\le T}$ is non-increasing and $K_{0}$ is bounded.

\end{Remark}
Note that our standing assumptions \reff{eq: ass f}-\reff{eq: ass g}-\reff{eq: ass K int 0} ensure that  \BSDE$_{K}(f,g,\zeta)$ admits a trivial super-solution
\be\label{eq: trivial supersol}
({\rm y}_{t},{\rm z}_{t})=((1+(T-t))L,0),
\ee
which is bounded.
In particular, \cite[Theorem 4.2]{peng1999monotonic} implies that  \BSDE$_{K}(f,g,\zeta)$ admits a minimal super-solution. We  denote it by $(\Ybb^{\zeta},\Zbb^{\zeta})$, and let  $\Kbb^{\zeta}$ be  the non-decreasing process defined on $[0,T]$ by
$$
\Kbb^{\zeta}_{\tau_{\zeta}}=0 \mbox{ and } \Ybb^{\zeta}_{\cdot \vee \tau_{\zeta}}=g(X^{\zeta}_{T})+\int_{\cdot\vee \tau_{\zeta}}^{T} f_{s}(X^{\zeta}_{s},\Ybb^{\zeta}_{s},\Zbb^{\zeta}_{s})ds - \int_{\cdot\vee \tau_{\zeta}}^{T} \Zbb^{\zeta}_{s} dW_{s}+\Kbb^{\zeta}_{T}-\Kbb^{\zeta}_{\cdot \vee \tau_{\zeta}}.
$$

We do not impose any Lipschitz continuity assumption on $g$, although it is used in the  unconstrained case $K=\R^{d}$ to obtain the Lipschitz continuity of the map $\xi\mapsto \hat \Yc^{(\tau,\xi)}_{\tau}$.
 Instead, we  assume that the map
\be\label{eq: ass hat g}
\hat g : x\in \R^{d}\mapsto \sup_{u\in \R^{d}}\left(g(x+u)-\delta_{T}(u)\right), \; x\in \R^{d}
\mbox{ is $L$-Lipschitz continuous.}
\ee
This map is usually referred to as the 'face-lift' of $g$ for the constraint  $K_{T}$, compare with e.g.~\cite{Broadie,bouchard2010obstacle,cvitanic1999super}. We shall see below that it provides the correct time-$T$ boundary condition for our constrained backward differential equation.
Intuitively, this means that assuming that $g$ is Lipschitz is useless, whenever $\hat g$ is, which is a weaker condition
\footnote{The fact  that $\hat g$ inherits the Lipschitz--continuity property from $g$ is by construction,
whereas the converse is not valid: for $d=1$,  $K=\R_{+}$ and  {$g:x\in\R_+\longmapsto \mathbf1_{\{x< 1\}}$},  we have  $\hat g:x\in\R_+\longmapsto 1$.}.
\vs2

Our  main result  shows that the map $\zeta\in \Db_{2}\mapsto \hat \Yc^{\zeta}_{\tau_{\zeta}}$ satisfies similar regularity properties in time and space as in the unconstrained case. It also shows that the non-decreasing process $\hat \Kc^{\zeta}$ is continuous on $[0,T)$ with a final jump of size $(\hat g-g)(X^{\zeta}_{T})$. In particular, $\hat \Yc^{\zeta}_{T-}=\hat g(X^{\zeta}_{T})$ on $\{\tau_{\zeta}<T\}$.

From now on, we denote by $C_{L}$ a generic constant which depends only on $L$, and may change from line to line.

\begin{theorem}\label{thm: regularity constrained bsde} The following holds for  all  $\zeta,\zeta'\in \Db_{2}$:
\begin{enumerate}[{\rm (a)}]
\item If $\tau_{\zeta}\le \tau_{\zeta'} {<T}$, then
\be\label{eq: thm main bade lip space holder time}
|\Ybb^{\zeta}_{\tau_{\zeta}}-\E_{\tau_{\zeta}}[\Ybb^{\zeta'}_{\tau_{\zeta'}}]|  \le C_{L} \left(\E_{\tau_{\zeta}}[|\tau_{\zeta'}-\tau_{\zeta}|]^{\frac12}+\E_{\tau_{\zeta}}[|\xi_{\zeta'}-\xi_{\zeta}|]\right).
\ee
\item If $\tau:=\tau_{\zeta}=\tau_{\zeta'} {<T}$, then
\be\label{eq: thm main bade delta lip space}
-\delta_{\tau}(\xi_{\zeta'}-\xi_{\zeta})\le \Ybb^{\zeta}_{\tau}-\Ybb^{\zeta'}_{\tau}  \le \delta_{\tau}(\xi_{\zeta}-\xi_{\zeta'}).
\ee
\item If $\tau_{\zeta}<T$,  then  $\Kbb^{\zeta}_{\cdot \wedge \vartheta}$ has  continuous path  for each stopping time $\vartheta<T$. Moreover, if $(\vartheta_{n})_{n\ge 1}$ is a   sequence of stopping times with values in $[\tau_{\zeta},T)$ such that $\vartheta_{n}\to T$, then
$
 \Ybb^{\zeta}_{\vartheta_{n}}\to \hat g(X^{\zeta}_{T}).
$
\end{enumerate}
\end{theorem}

Sections \ref{sec: estimates via the strong dual formulation} and \ref{sec: weak versus strong} are devoted to the proof of these results.  In view of Proposition \ref{prop: equiv strong and weak}  and Theorem \ref{thm: Doob Meyer decomposition}: {\rm (a)}~is a consequence of Corollary \ref{cor: space regularity} and Proposition \ref{prop: stab time}, {\rm (b)}~follows from Proposition \ref{prop: Yb self face lifted},   and Proposition \ref{prop: path continuity} implies {\rm (c)}.
\section{Estimates via the strong dual formulation}\label{sec: estimates via the strong dual formulation}

As explained in the introduction, the constrained backward differential equation \BSDE$_{K}(f,$ $g,$ $\zeta)$ admits a dual    representation which is formulated as an optimal control problem on a family of unconstrained backward stochastic differential equations written under a suitable family of equivalent laws, see Section \ref{sec: weak versus strong} for a precise formulation. Although, each unconstrained backward stochastic differential equation satisfies the usual H\"older and Lipschitz regularity properties, this does not seem to allow one to obtain the estimates of Theorem \ref{thm: regularity constrained bsde}. The reason is that the optimal control problem is of singular type: constants may blow up when passing to the supremum.

The main idea of this paper is to start with a strong version of this dual optimal control problem. Strong meaning that the probability measure is fixed, but we incorporate the control directly in the dynamics. It turns out to be much more flexible. In particular, space stability is essentially trivial in this setting, see Corollay \ref{cor: space regularity}. More importantly, we can show that the corresponding value process is itself automatically 'face-lifted', see Proposition \ref{prop: Yb self face lifted} below. This will be the key result to obtain the time regularity estimates of Proposition \ref{prop: stab time}.

\subsection{The strong dual formulation}

Let $\Uc$ denote the collection of $\R^{d}$-valued bounded predictable processes. Note that  $\delta(\nu)$ is bounded for each $\nu \in \Uc$, see \reff{eq: ass K borne}. To each $(\zeta,\nu)\in \Db_{2}\x \Uc$, we associate the stochastic driver
$$
f^{\zeta,\nu}:(t,y,z)\in [0,T]\x \R\x \R^{d}\mapsto  \left(f_{t}(X^{\zeta,\nu}_{t},y,z)-\delta_{t}(\nu_{t})\right)\1_{[\![\tau_{\zeta},T]\!]}(t)
$$
where  $X^{\zeta,\nu}$ is the solution of
\be
X^{\zeta,\nu}&=&\xi_{\zeta}+\int_{\tau_{\zeta}}^{\cdot\vee \tau_{\zeta}} \left(b_{s}(X^{\zeta,\nu}_{s})+\nu_{s}\right)ds +  \int_{\tau_{\zeta}}^{\cdot\vee \tau_{\zeta}} \sigma_{s}(X^{\zeta,\nu}_{s})dW_{s}.
\label{eq: def bar X}
\ee
Given  $ \tau \in \Tc_{\tau_{\zeta}}, \vartheta \in \Tc_{\tau}$ and $G\in \Lb_{2}(\Fc_{\vartheta})$, we set
$$
\Ec^{\zeta,\nu}_{\tau,\vartheta}[G]:=U_{\tau}
$$
where $(U,V)\in \Sb_{2}\x \Hb_{2}$ is the solution of
\be\label{eq: BSDE f zeta nu for non line expect}
U &=& G+\int_{\cdot\vee\tau_{\zeta}}^{\vartheta} f^{\zeta,\nu}_{s}(U_{s},V_{s}) ds-\int_{\cdot\vee\tau_{\zeta}}^{\vartheta}  V_{s}dW_{s}\; \mbox{ on } [0,T].
\ee
In the special case where $\vartheta\equiv T$ and $G=g(X^{\zeta,\nu}_{T})$, the solution of \reff{eq: BSDE f zeta nu for non line expect} is denoted by $(Y^{\zeta,\nu},Z^{\zeta,\nu})$. In particular,
 $$
 Y^{\zeta, \nu}_{\cdot}=\Ec^{\zeta,\nu}_{\cdot,T}[g(X^{\zeta,\nu}_{T})].
 $$
We next define our optimal control problem
\be\label{eq: def Yc}
\Yc^{\zeta}_{\tau}= \esssup \{ Y^{\zeta,\nu}_{\tau},\nu \in \Uc,\;\nu\1_{\[0,\tau\]}\equiv 0\}\;,\; \zeta \in \Db_{2},\;\tau \in \Tc_{\tau_{\zeta}}.
\ee
 Note that $\Yc^{\zeta}_{\tau}=\Yc^{(\tau,X^{\zeta}_{\tau})}_{\tau}$ since $X^{\zeta,0}=X^{\zeta}$.

\begin{remark}\label{rem: bound bar Y et Yc} The conditions \reff{eq: ass b sigma}-\reff{eq: ass f}-\reff{eq: ass g} imply that $Y^{\zeta,\nu}$ is bounded in $\Lb_{\infty}$  uniformly in $\zeta\in \Db_{2}$,  for all $\nu \in \Uc$,  see Lemma \ref{lem: appendix linearisation}. Moreover, $Y^{\zeta,\nu}\le {\rm y}$ defined in  \reff{eq: trivial supersol}, for all $\nu \in \Uc$, see Lemma \ref{lem: appendix stabilite et comp ?}.
In particular,
$Y^{\zeta, 0}\le \Yc^{\zeta}\le {\rm y}$, so that $\Yc^{\zeta}$ is bounded in $\Lb_{\infty}$ uniformly in $\zeta\in \Db_{2}$. The bound depends only on $L$.
\end{remark}
\subsection{Terminal face-lift and space stability}

The first result of this section concerns the face-lift of the terminal condition. It shows that $g$ can be replaced by $\hat g$ in  \reff{eq: def Yc}. Apart from being of self-interest, this property will be used latter on in the proof of the space stability in our setting where $\hat g$ is assumed to be Lipschitz while $g$ may not be. It will also be used to characterize the limit $\lim_{t\uparrow T} \Yc^{\zeta}_{t}$.

\begin{Proposition}\label{prop: facelift} $\Yc^{\zeta}_{\tau_{\zeta}}= \esssup \{ \Ec^{\zeta,\nu}_{\tau_{\zeta},T}[\hat g(X^{\zeta,\nu}_{T})],\nu \in \Uc\}$  on $\{\tau_{\zeta}<T\}$, for all $\zeta \in \Db_{2}$.
\end{Proposition}

\proof Since $\hat g\ge g$ by construction, one inequality is trivially deduced from Lemma \ref{lem: appendix stabilite et comp ?}.
We therefore concentrate on the difficult inequality. Fix $\zeta:=(\tau,\xi)\in \Db_{2}$. For sake of simplicity, we assume that $\tau<T$ a.s., the general case is handled by using the fact that $\{\tau<T\}=\cup_{n\ge 1} \{\tau\le T-n^{-1}\}$.    For $\nu \in \Uc$, and $u\in  {\Lb_{\infty}(\Fc_{T-\eps_{\circ}})}$, for some $\eps_{\circ}>0$, we then define
$$
\tau_{\eps}:=(T-\eps)\vee \tau\;,\;    \nu^{\eps}:=\nu\1_{\[\tau,\tau_{\eps}\]}+\frac{u}{T-\tau_{\eps}}\1_{\]\tau_{\eps},T\]},\; 0<\eps<\eps_{\circ}.
$$
Then,  \reff{eq: def Yc} combined with the tower property for non linear expectations imply that
\b*
  \Yc^{\zeta}_{\tau }&\ge&\Ec_{\tau,\tau_{\eps}}^{\zeta,\nu}\left[  Y^{\zeta,\nu^{\eps}}_{\tau_{\eps}}\right].
\e*
We claim that, after possibly considering a subsequence,
\be\label{eq: claim prop: facelift}
\liminf_{\eps\to 0}\Ec_{\tau,\tau_{\eps}}^{\zeta,\nu}\left[  Y^{\zeta,\nu^{\eps}}_{\tau_{\eps}}\right]\ge \Ec_{\tau,T}^{\zeta,\nu}\left[   g (X^{\zeta,\nu}_{T}+u)-\delta_{T}(u)  \right].
\ee
By arbitrariness of $\eps_{\circ}>0$ and $u\in  {\Lb_{\infty}(\Fc_{T-{\eps_{\circ}}})}$, this implies that
\b*
  \Yc^{\zeta}_{\tau }&\ge& \esssup_{u\in \Lb_{\infty}(\Fc_{T-})} \Ec_{\tau,T}^{\zeta,\nu}\left[   g (X^{\zeta,\nu}_{T}+u)-\delta_{T}(u)  \right].
\e*
Since  the map $(x,u)\in \R^{d}\x \R^{d} \mapsto  g (x+u)-\delta_{T}(u)$ is Borel, it follows from  \cite[Proposition 7.40, p184]{BeSh78}, that, for each $\iota>0$, we can find a   universally measurable map $x\in \R^{d}\mapsto \tilde  u_{\iota}(x)$ such that $\hat g(x)$ $=$ $\sup\{  g (x+u)-\delta_{T}(u),$ $u\in \R^{d}\}\le g (x+\tilde u_{\iota}(x))-\delta_{T}(\tilde u_{\iota}(x)) +\iota$ for all $x\in \R^{d}$. By  \cite[Lemma 7.27, p173]{BeSh78}, we can find a Borel measurable map $x\in \R^{d}\mapsto \hat  u_{\iota}(x)$ such that $\hat u_{\iota}(X^{\zeta,\nu}_{T})=\tilde u_{\iota}(X^{\zeta,\nu}_{T})$ a.s.  Since $X^{\zeta,\nu}_{T}\in \Lb_{2}(\Fc_{T-})$ by left-continuity of its path, this implies that
\b*
  \Yc^{\zeta}_{\tau }&\ge& \Ec_{\tau,T}^{\zeta,\nu}\left[  \hat g (X^{\zeta,\nu}_{T})\1_{\{|\hat u_{\iota}(X^{\zeta,\nu}_{T})|\le n\}}-\iota +   g (X^{\zeta,\nu}_{T})\1_{\{|\hat u_{\iota}(X^{\zeta,\nu}_{T})|> n\}}  \right],
\e*
for each $n\ge 1$ and $\iota \in (0,1)$.
The required result then follows from  the stability principle, see Lemma \ref{lem: appendix stabilite et comp ?},  by sending first $n\to \infty$ and then $\iota\to 0$.

It remains to prove our claim \reff{eq: claim prop: facelift}.  We first deduce from Lemma \ref{lem : X+u sur temps eps} stated at the end of this section that, after possibly passing to a subsequence,
\b*
X^{\zeta,\nu^{\eps}}_{T} \to X^{\zeta,\nu}_{T} +u \;\;a.s.\; \mbox{ as  } \eps\to 0,
\e*
and therefore
\b*
\liminf_{\eps\to 0} g(X^{\zeta,\nu^{\eps}}_{T}) \ge g( X^{\zeta,\nu}_{T} +u),
\e*
by lower-semicontinuity of $g$.
Moreover, if  $(H^{\eps})_{\eps>0}$ is a sequence of positive   processes such that
$\sup_{[\tau_{\eps},T]}H^{\eps}\to 1$ a.s.~as $\eps\to 0$, then
\b*
\liminf_{\eps \to 0} \int_{\tau_{\eps}}^{T} H^{\eps}_{s} \left(  f_{s}(X^{\zeta,\nu^{\eps}}_{s},0)-\delta_{s}(\nu^{\eps}_{s})\right)ds
 &\ge&  -
 \lim_{\eps \to 0}  \; (T-\tau_{\eps})^{-1}\int_{\tau_{\eps}}^{T}H^{\eps}_{s} \delta_{s}(u)ds\\
  &\ge&  -\delta_{T}(u)
\e*
since   $t\mapsto \delta_{t}(u)$ is left-continuous at $T$ and $f(\cdot,0)$ is bounded.
We can then  combine Lemma \ref{lem: appendix linearisation} and Lemma \ref{lem: appendix liminf eps and E condieps}  to obtain
$$
\liminf_{\eps\to 0}  Y^{\zeta,\nu^{\eps}}_{\tau_{\eps}}\ge g(X^{\zeta,\nu}_{T} +u)-\delta_{T}(u)=:G,
$$
after possibly passing to a subsequence.  
Now observe that $\Ec_{\cdot,\tau_{\eps}}^{\zeta,\nu}[Y^{\zeta,\nu^{\eps}}_{\tau_{\eps}}]$ coincides with the first component of the backward differential equation which driver is
$f^{\zeta,\nu}\mathbf1_{[\![0,\tau_\eps]\!]}$ and which terminal condition  is $Y^{\zeta,\nu^{\eps}}_{\tau_{\eps}}$ at $T$.  {It is constant on $[\![\tau_\eps,T]\!]$.}
Then, by
 the stability  and comparison principles in Lemma \ref{lem: appendix stabilite et comp ?},   we  obtain
\b*
\Ec_{\tau,\tau_{\eps}}^{\zeta,\nu}\left[  Y^{\zeta,\nu^{\eps}}_{\tau_{\eps}}\right]
&\ge& \Ec_{\tau,T}^{\zeta,\nu}\left[  Y^{\zeta,\nu^{\eps}}_{\tau_{\eps}}\right] -C_{1} \E_{\tau}\left[\int_{\tau_{\eps}}^{T} |f^{\zeta,\nu}_{s}(Y^{\zeta,\nu^{\eps}}_{\tau_{\eps}},0)|^{2} ds\right]^{\frac12}
\\
&\ge&  \Ec_{\tau,T}^{\zeta,\nu}\left[  G\right] -C_{2}\left( \E_{\tau}\left[\int_{\tau_{\eps}}^{T} |f^{\zeta,\nu}_{s}(Y^{\zeta,\nu^{\eps}}_{\tau_{\eps}},0)|^{2} ds\right]^{\frac12} + \E_{\tau}\left[\{ (Y^{\zeta,\nu^{\eps}}_{\tau_{\eps}}-G)^{-}\}^{2}\right]^{\frac12} \right),
\e*
in which the constants $C_{1}$ and $C_{2}$ do not depend on $(u,\eps)$,  {see Remark \ref{rem: bound bar Y et Yc}}.
Since $\nu \in \Uc$, it follows from Lemma \ref{lem: appendix linearisation} {together with \reff{eq: ass b sigma} and \reff{eq: ass f}} that $(Y^{\zeta,\nu^{\eps}}_{\tau_{\eps}},f^{\zeta,\nu}(Y^{\zeta,\nu^{\eps}}_{\tau_{\eps}},0)\1_{[\![\tau_{\eps},T]\!]})$ is bounded in $\Lb_{2}\times \Sb_{2}$ uniformly in $\eps>0$. Then, combining the above shows that, {along a subsequence if necessary}, the right-hand side term in the last inequality converges to $0$ {as $\eps\downarrow0$}.
{Hence, the claim \reff{eq: claim prop: facelift}
holds, which completes the proof.}
\ep
\\

Since $\hat g$ is assumed to be Lipschitz continuous, the  stability in space is now an easy consequence of the representation given in Proposition \ref{prop: facelift}
\begin{Corollary}\label{cor: space regularity}
$
\left| \Yc^{\zeta_{1}}_{\tau}- \Yc^{\zeta_{2}}_{\tau}\right|\le  C_{L} |\xi_{\zeta_{1}}-\xi_{\zeta_{2}}|,
$
for all $\zeta_{1},\zeta_{2}\in \Db_{2}$  with $\tau:=\tau_{\zeta_{1}}=\tau_{\zeta_{2}} {<T}$.
\end{Corollary}
\proof We simply use the fact that
$$
\left| \Yc^{\zeta_{1}}_{\tau}- \Yc^{\zeta_{2}}_{\tau}\right|\le \esssup_{\nu \in \Uc}\left|\Ec^{\zeta_{1},\nu}_{\tau,T}[\hat g(X^{\zeta_{1},\nu}_{T})]-\Ec^{\zeta_{2},\nu}_{\tau,T}[\hat g(X^{\zeta_{2},\nu}_{T})]\right|
$$
by Proposition \ref{prop: facelift}.  The right-hand side is bounded by $  |\xi_{\zeta_{1}}-\xi_{\zeta_{2}}|$ up to a multiplicative constant under our Lipschitz continuity assumptions  \reff{eq: ass b sigma}-\reff{eq: ass f}-\reff{eq: ass hat g}, see Lemma \ref{lem: appendix stabilite et comp ?}.
\ep
\\

We conclude this section with the technical lemma that was used  in the {proof of Proposition \ref{prop: facelift}}. 
The proof is trivial under \reff{eq: ass b sigma} and we  omit it.  It is not difficult to see that it remains correct without the boundedness assumption on $(b,\sigma)$, they only need to be Lipschitz continuous in space, uniformly in time (but then the constant appearing in the bound  depends on $(u,\zeta)$ as well).

\begin{Lemma}\label{lem : X+u sur temps eps} Fix $\zeta\in \Db_{2}$,  $\vartheta\in \Tc_{\tau_{\zeta}}$  and  $u\in \Lb_{\infty}(\R^{d},\Fc_{\tau_{\zeta}})$. Set
$
 \nu:=\eps^{-1} u  \1_{\{\eps>0\}}\1_{[\![\tau_{\zeta},\vartheta]\!]} ,
$
with $\eps:=\vartheta-\tau_{\zeta}$.
Then,
$$
\sup_{t\le T}\E_{\tau_{\zeta}}\left[\left| X^{\zeta,\nu}_{ t\vee \tau_{\zeta} \wedge \vartheta}-\xi_{\zeta}-\eps^{-1}u(t\vee \tau_{\zeta}\wedge \vartheta-\tau_{\zeta})\1_{\{\eps>0\}}\right|^{2}\right]\le C_{L}\E_{\tau_{\zeta}}[\eps].
$$
\end{Lemma}

\subsection{Dynamic programming and face-lifting on $[0,T)$}

We first recall the dynamic programming principle for the optimal control problem \reff{eq: def Yc}. It will be used later on in this section to prove that the value process is automatically face-lifted,  {see} Proposition \ref{prop: Yb self face lifted}.

\begin{Proposition}\label{prop: dpp sup bar Y}  For all $\zeta\in \Db_{2}$ and $\vartheta \in \Tc_{\tau_{\zeta}}$,
 $$
 \Yc^{\zeta}_{\tau_{\zeta}}= \esssup_{\nu \in \Uc} \Ec^{\zeta,\nu}_{\tau_{\zeta},\vartheta}[\Yc^{(\vartheta,X^{\zeta,\nu}_{\vartheta})}_{\vartheta}].
 $$
\end{Proposition}
\proof The proof is standard. Since
$
Y^{\zeta,\nu}_{\tau_{\zeta}}=\Ec^{\zeta,\nu}_{\tau_{\zeta}}[Y^{(\vartheta,X^{\zeta,\nu}_{\vartheta}),\nu}_{\vartheta}]\le \Ec^{\zeta,\nu}_{\tau_{\zeta}}[\Yc^{(\vartheta,X^{\zeta,\nu}_{\vartheta})}_{\vartheta}]
$
by the tower property for non-linear expectations and by the comparison principle, see Lemma \ref{lem: appendix stabilite et comp ?}, one inequality is trivial. As for the reverse inequality, we observe  that the family
$\{Y^{(\vartheta,X^{\zeta,\nu}_{\vartheta}),\nu'}_{\vartheta},\nu' \in \Uc\}$ is directed upward. Then \cite[Proposition VI.1.1]{Ne75} ensures that we can find a sequence $(\nu'_{n})_{n\ge 1}\subset \Uc$ such that
$Y^{n}_{\vartheta}:=Y^{(\vartheta,X^{\zeta,\nu}_{\vartheta}),\nu'_{n}}_{\vartheta}\uparrow \Yc^{(\vartheta,X^{\zeta,\nu}_{\vartheta})}_{\vartheta}=:\Yc_{\vartheta}$ a.s.~as $n\to \infty$. Since,
$Y^{1}_{\vartheta}$ and $\Yc_{\vartheta}$ are bounded in $\Lb_{2}$, see Remark \ref{rem: bound bar Y et Yc},  the convergence holds in $\Lb_{2}$ as well. Moreover, we can find a constant $C>0$ which does not depend on $n$ and such that
\b*
 \Ec^{\zeta,\nu}_{\tau_{\zeta},\vartheta}[Y^{n}_{\vartheta}]
 &\ge&
 \Ec^{\zeta,\nu}_{\tau_{\zeta},\vartheta}[\Yc_{\vartheta}]-C\; \E_{\tau_{\zeta}}[|Y^{n}_{\vartheta}-\Yc_{\vartheta}|^{2}]^{\frac12},
\e*
see Lemma \ref{lem: appendix stabilite et comp ?}.   {The latter combined with \reff{eq: def Yc}
implies} that
\b*
 \Yc^{\zeta}_{\tau_{\zeta}}\ge  \Ec^{\zeta,\nu}_{\tau_{\zeta},\vartheta}[\Yc_{\vartheta}]- \lim_{n\to \infty}C\;\E_{\tau_{\zeta}}[|Y^{n}_{\vartheta}-\Yc_{\vartheta}|^{2}]^{\frac12}
= \Ec^{\zeta,\nu}_{\tau_{\zeta},\vartheta}[\Yc_{\vartheta}],
\e*
and we conclude by arbitrariness of $\nu \in \Uc$.
\ep
\\

We can now show that $\xi\in \Lb_{2}(\Fc_{\tau})\mapsto \Yc^{(\tau,\xi)}_{\tau}$ is itself automatically face-lifted  in the following sense.

 \begin{Proposition}\label{prop: Yb self face lifted}
  For all $\zeta\in \Db_{2}$,
 \be\label{eq: face lift t<T}
 \Yc^{\zeta}_{\tau_{\zeta}}= \esssup_{u\in \Lb_{\infty}(\R^{d},\Fc_{\tau_{\zeta}})}(\Yc^{(\tau_{\zeta},{\xi_{\zeta}+u})}_{\tau_{\zeta}}-\delta_{{\tau_{\zeta}}}(u))\;\mbox{ a.s. on $\{\tau_{\zeta}<T\}$.}
 \ee
 \end{Proposition}

 \proof Take $\zeta:=(\tau,\xi)\in \Db_{2}$. One inequality follows from the fact that $\delta_{\tau}(0)=0$. Fix $u\in \Lb_{\infty}(\Fc_{\tau})$, $\eps>0$, and set $\tau_{\eps}:=(\tau+\eps)\wedge T$ and $\nu^{\eps}:= \eps^{-1} u\1_{\[\tau,\tau_\eps\]}$. It follows from Lemma \ref{lem : X+u sur temps eps}  that
 $$
 \E_{\tau}[|X^{\zeta,\nu^\eps}_{\tau_\eps} - \xi - u|^{2}]^{\frac12}\le C_{L}\; \eps^{\frac12}.
 $$
 Then, by appealing to Proposition \ref{prop: dpp sup bar Y}, Lemma \ref{lem: appendix linearisation}, \reff{eq: ass f}  and  Corollary \ref{cor: space regularity} successively, we can find a family of non-negative continuous processes $(H^{\eps})_{\eps >0}$, uniformly bounded in $\Sb_{2}$, such that  $H^{\eps}_{\tau_\eps}\to 1$ in $\Lb_{1}$ as $\eps\to 0$ and
\b*
  \Yc^{\zeta}_{\tau}&\ge&\E_\tau\left[H^{\eps}_{\tau_\eps} \Yc^{(\tau_\eps,X^{\zeta,\nu^{\eps}}_{\tau_\eps})}_{\tau_\eps} +\int_{\tau}^{\tau_\eps} H^{\eps}_{s} f^{\zeta,\nu^{\eps}}_{s}(0)ds\right]
\\
&\ge&
\E_\tau\left[H^{\eps}_{\tau_\eps} \Yc^{(\tau_\eps,X^{\zeta,\nu^{\eps}}_{\tau_\eps})}_{\tau_\eps} -\int_{\tau}^{\tau_\eps} H^{\eps}_{s} {\delta_s(\nu^{\eps}_{s})}ds\right]-C_{L}\eps
\\
&\ge&
\E_\tau\left[H^{\eps}_{\tau_\eps} \Yc^{(\tau_\eps,\xi+u)}_{\tau_\eps}- \eps^{-1}\int_{\tau}^{\tau_\eps} H^{\eps}_{s} \delta_s(u)ds \right]-C_{L} \eps^{\frac12}.
\e*
Since $H^{\eps}\ge 0$, we can use \reff{eq : ass delta non-increasing} and Remark \ref{rem: bound bar Y et Yc} to obtain
\b*
  \Yc^{\zeta}_{\tau} &\ge&
\E_\tau\left[H^{\eps}_{\tau_\eps} \Yc^{(\tau_\eps,\xi+u)}_{\tau_\eps}-\delta_\tau(u) \eps^{-1}\int_{\tau}^{\tau_\eps} H^{\eps}_{s} ds \right]-C_{L} \eps^{\frac12}\\
&\ge& \E_\tau\left[\Yc^{(\tau_\eps,\xi+u)}_{\tau_\eps}-\delta_\tau(u) \eps^{-1}\int_{\tau}^{\tau_\eps} H^{\eps}_{s} ds \right]-C_{L} (\eps^{\frac12}+\E_{\tau}[|H^{\eps}_{\tau_\eps} -1|])
\e*
where $\eps^{-1}\int_{\tau}^{\tau_\eps} H^{\eps}_{s} ds$ and $H^{\eps}_{\tau_{\eps}}$ converge to $1$ in $\Lb_{1}$, so that
\b*
  \Yc^{\zeta}_{\tau} &\ge&\liminf_{\eps\to 0}
\E_\tau\left[  \Yc^{(\tau_\eps,\xi+u)}_{\tau_\eps} \right]-\delta_\tau(u).
\e*
It remains to show that
\beq\label{eq: liminf intermediary}
\liminf_{\eps\to 0} \E_\tau\left[   \Yc^{(\tau_\eps,\xi+u)}_{\tau_\eps}\right]\ge  \Ec^{(\tau ,\xi+u),\nu}_{\tau,T}[\hat g(X^{(\tau ,\xi+u),\nu}_{T})].
\eeq
Then, the arbitrariness of $\nu\in \Uc$   will allow one to conclude by  appealing to Proposition \ref{prop: facelift}.
In order to alleviate notations, we  set $\hat Y^{\tau_{1},\xi_{1}}_{\tau_{2}}:=\Ec^{(\tau_{1},\xi_{1}),\nu}_{\tau_{2},T}[\hat g(X^{(\tau_{1},\xi_{1}),\nu}_{T})]$
for any $\tau_{1}, \tau_{2}\in \Tc$ and $\xi_{1} \in \Lb_{2}(\Fc_{\tau_{1}})$. By Proposition \ref{prop: facelift},
\b*
\E_\tau[   \Yc^{(\tau_\eps,\xi+u)}_{\tau_\eps}]\ge \E_\tau[  \hat Y_{\tau_{\eps}}^{\tau_{\eps},\xi+u}]=  \hat Y^{\tau,\xi+u}_{\tau} + \E_\tau[\hat Y^{\tau_\eps,\xi+u}_{\tau_\eps}] -\hat Y^{\tau,\xi+u}_{\tau}.
\e*
We now observe that
\b*
\E_\tau[\hat Y^{\tau_\eps,\xi+u}_{\tau_\eps}] -\hat Y^{\tau,\xi+u}_{\tau}
&=&  \E_\tau[\hat Y^{\tau_{\eps},\xi+u}_{\tau_{\eps}}-\hat Y^{\tau_\eps,X^{(\tau,\xi+u),\nu}_{\tau_{\eps}}}_{\tau_\eps} ]\\
&& +\;\E_\tau[\hat Y^{\tau_\eps,X^{(\tau,\xi+u),\nu}_{\tau_{\eps}}}_{\tau_\eps}]-\Ec^{(\tau,\xi+u),\nu}_{\tau,\tau_{\eps}}[\hat Y^{\tau_\eps,X^{(\tau,\xi+u),\nu}_{\tau_{\eps}}}_{\tau_\eps}]
\e*
in which, by Lemma \ref{lem: appendix stabilite et comp ?} combined with \reff{eq: ass b sigma}-\reff{eq: ass f}-\reff{eq: ass hat g},
$$
\lim_{\eps\to 0} \E_\tau[\hat Y^{\tau_{\eps},\xi+u}_{\tau_{\eps}}-\hat Y^{\tau_\eps,X^{(\tau,\xi+u),\nu}_{\tau_{\eps}}}_{\tau_\eps} ]=0.
$$
By the same assumptions combined with Lemma \ref{lem: appendix linearisation},
 $$
 \lim_{\eps\to 0} \E_\tau[\hat Y^{\tau_\eps,X^{(\tau,\xi+u),\nu}_{\tau_{\eps}}}_{\tau_\eps}]-\Ec^{(\tau,\xi+u),\nu}_{\tau,\tau_{\eps}}[\hat Y^{\tau_\eps,X^{(\tau,\xi+u),\nu}_{\tau_{\eps}}}_{\tau_\eps}]=0.
 $$
 This proves \reff{eq: liminf intermediary} and completes the proof.
 \ep

\subsection{Stability in time}

We now turn to the proof of the stability in time. The lower estimate trivially follows from the dynamic programming principle of Proposition \ref{prop: dpp sup bar Y}. The second one is much more delicate. It is obtained by
a suitable use of the face-lifting phenomenon  observed in Proposition \ref{prop: Yb self face lifted}. It allows one to absorb the singularity due to the control when passing to the supremum over $\Uc$, see \reff{eq: use face lift for time regu} below.

\begin{Proposition}\label{prop: stab time} For all $\zeta \in \Db_{2}$ and  $\vartheta \in \Tc_{\tau_{\zeta}}$
$$
\left|\Yc^{\zeta}_{\tau_{\zeta}}-\E_{\tau_{\zeta}}[\Yc^{(\vartheta,\xi_{\zeta})}_{\vartheta}]\right|\le C_{L}\;\E_{\tau_{\zeta}}[\vartheta-\tau_{\zeta}]^{\frac12}.
$$

\end{Proposition}

\proof 1.  By Proposition \ref{prop: dpp sup bar Y}, Remark \ref{rem: bound bar Y et Yc}, \reff{eq: ass f} and Lemma \ref{lem: appendix linearisation}
$$
\Yc^{\zeta}_{\tau_{\zeta}}\ge \Ec_{\tau_\zeta,\vartheta}^{\zeta,0}[ \Yc^{(\vartheta,X^{\zeta,0}_{\vartheta})}_{\vartheta}]\ge \E_{\tau_{\zeta}}[\Yc^{(\vartheta,X^{\zeta,0}_{\vartheta})}_{\vartheta}]-C_{L}\E_{\tau_{\zeta}}[\vartheta-\tau_{\zeta}]^{\frac12},
$$
while Corollary \ref{cor: space regularity} and \reff{eq: ass b sigma} imply
$$
\E_{\tau_{\zeta}}[\Yc^{(\vartheta,X^{\zeta,0}_{\vartheta})}_{\vartheta}]\ge \E_{\tau_{\zeta}}[\Yc^{(\vartheta,\xi_{\zeta})}_{\vartheta}]-C_{L}\E_{\tau_{\zeta}}[\vartheta-\tau_{\zeta}]^{\frac12}.
$$

2. We now turn to the reverse inequality. Set $\zeta=(\tau,\xi)$ and let
$U^{\zeta,\nu} :=\beta Y^{\zeta,\nu}_{\cdot\vee \tau}-\int_{\tau}^{\cdot\vee \tau} \beta_{s}\delta_{s}(\nu_{s})ds$, 
where $\beta_{s}:=e^{L (s-\tau)^{+}}$, recall \reff{eq: ass f}. Then,
\b*
  U^{\zeta,\nu}_{t_{1}}&{=} &   U^{\zeta,\nu}_{t_{2}} +\int_{\tau\vee t_{1}}^{\tau\vee t_{2}}\beta_{s}\left(f_s(X^{\zeta,\nu}_{s},Y^{\zeta,\nu}_{s},Z^{\zeta,\nu}_{s}) -L Y^{\zeta,\nu}_{s}\right) ds-\int_{\tau\vee t_{1}}^{\tau\vee t_{2}} \beta_{s}Z^{\zeta,\nu}_{s}dW_{s}
\e*
for $t_{2}\ge t_{1}$.
Since  {$y\in \R\mapsto f(\cdot,y,\cdot)-L y$ is non-increasing} by \reff{eq: ass f}, and  $\delta\ge 0$ by \reff{eq: ass K int 0}, it follows that
\beq\label{eq: two last}
  U^{\zeta,\nu}_{t_{1}}\le   U^{\zeta,\nu}_{t_{2}} +\int_{\tau\vee t_{1}}^{\tau\vee t_{2}}\left(\beta_{s} f_{s}(X^{\zeta,\nu}_{s},  \beta_{s}^{-1} U^{\zeta,\nu}_{s},Z^{\zeta,\nu}_{s}) -L U^{\zeta,\nu}_{s} \right)ds-\int_{\tau\vee t_{1}}^{\tau\vee t_{2}}\beta_{s}Z^{\zeta,\nu}_{s}dW_{s}.
\eeq
 On the other hand, we  can use \reff{eq: def Yc}, the fact that $\delta\ge 0$ is sublinear and $\beta \ge 1$, \reff{eq : ass delta non-increasing}, Remark \ref{rem: bound bar Y et Yc} and Proposition \ref{prop: Yb self face lifted} to deduce that
\be
U^{\zeta,\nu}_{\vartheta}&\le& \beta_{\vartheta} \Yc^{(\vartheta,X^{\zeta,\nu}_{\vartheta})}_{\vartheta}-\int_{\tau}^{\vartheta}\beta_{s} \delta_{s}(\nu_{s}) ds\nonumber \\
&\le& C_{L}(\beta_{\vartheta}-1)
+\Yc^{(\vartheta,X^{\zeta,\nu}_{\vartheta})}_{\vartheta}- \delta_{\vartheta}\left(\int_{\tau}^{\vartheta}\nu_{s}  ds\right)\nonumber\\
&\le&
C_{L}(\beta_{\vartheta}-1) + \Yc^{(\vartheta,X^{\zeta,\nu}_{\vartheta}-\int_{\tau}^{\vartheta} \nu_{s} ds)}_{\vartheta}.~~~\label{eq: use face lift for time regu}
\ee
The last inequality combined with \reff{eq: two last},
Lemma \ref{lem: appendix linearisation} and \reff{eq: ass f} leads to
\b*
Y^{\zeta,\nu}_{\tau}=  U^{\zeta,\nu}_{\tau}
&\le& \E_{\tau}\left[\hat H^{\tau}_{\vartheta} \left(C_{L}(\beta_{\vartheta}-1)+ \Yc^{(\vartheta,X^{\zeta,\nu}_{\vartheta}-\int_{\tau}^{\vartheta} \nu_{s} ds)}_{\vartheta}+C_{L}\;(\vartheta-\tau)\right)\right]
\e*
in which
$$
\hat H^{\tau}_{\vartheta}:=\sup_{[\tau,\vartheta]} e^{\int_{\tau}^{\cdot}(\kappa^{1}_{s}-2^{-1}|\kappa^{2}_{s}|^{2})ds +\int_{\tau}^{\cdot} \kappa_{s}^{2} dW_{s}}
$$
for some predictable processes $\kappa^{1}, \kappa^{2} $ that are uniformly bounded by a constant which only depends  on $L$.
We next use Corollary \ref{cor: space regularity}  together with  standard estimates, recall \reff{eq: ass b sigma},    to deduce from the above that
\b*
Y^{\zeta,\nu}_{\tau}- \E_{\tau}\left[\hat H^{\tau}_{\vartheta}   \Yc^{(\vartheta,\xi)}_{\vartheta}\right]
&\le& C_{L}\;\left(\E_{\tau}\left[\hat H^{\tau}_{\vartheta} \;|X^{\zeta,\nu}_{\vartheta}-\int_{\tau}^{\vartheta} \nu_{s} ds-\xi |\right]+\E_{\tau}\left[(\vartheta-\tau) \right]^{\frac12}\right)
\nonumber
\\
&\le&  C_{L}\;\E_{\tau}\left[(\vartheta-\tau) \right]^{\frac12}.\label{eq: proof upper bound time step 1}
\e*
Then, Remark \ref{rem: bound bar Y et Yc} implies that
$$
 \E_{\tau}\left[\hat H^{\tau}_{\vartheta}   \Yc^{(\vartheta,\xi)}_{\vartheta}\right]- \E_{\tau}\left[  \Yc^{(\vartheta,\xi)}_{\vartheta}\right]
 \le C_{L}\E_{\tau}\left[|\hat H^{\tau}_{\vartheta}-1|\right]\le C_{L}\E_{\tau}\left[(\vartheta-\tau) \right]^{\frac12}.
$$
Combining the two last inequalities and using the arbitrariness of $\nu \in \Uc$ leads to the required result.
\ep
\\

For later use, we state the following corollary of Proposition \ref{prop: stab time} and Corollary \ref{cor: space regularity}, recall \reff{eq: ass b sigma}.
\begin{Corollary}\label{cor: stab time on the path of X} For all $\zeta \in \Db_{2}$ and  $\vartheta \in \Tc_{\tau_{\zeta}}$
$$
\left|\Yc^{\zeta}_{\tau_{\zeta}}-\E_{\tau_{\zeta}}[\Yc^{(\vartheta,X^{\zeta}_{\vartheta})}_{\vartheta}]\right|\le C_{L}\;\E_{\tau_{\zeta}}[\vartheta-\tau_{\zeta}]^{\frac12}.
$$
\end{Corollary}

 \subsection{Path continuity and boundary limit}

The following is deduced from Proposition \ref{prop: facelift} and Corollary \ref{cor: stab time on the path of X}.

 \begin{Proposition}\label{prop: path continuity} Fix  $\zeta\in \Db_{2}$ such that $\tau_{\zeta}<T$. Then,  $\Yc^{\zeta}$   is continuous on $[0,T)$ and satisfies
  $$
 \Yc^{\zeta}_{T-}:=\lim_{s\uparrow T,s<T}  \Yc^{\zeta}_{s}=\hat g(X^{\zeta}_{T})\;.
 $$
  In particular, the process $\Yc^{\zeta}\1_{[0,T)}+\hat g(X^{\zeta}_{T})\1_{\{T\}}$ is continuous.
 \end{Proposition}

 \proof  1. We first show that $\Yc^{\zeta}$ is right-continuous.
 Fix $\vartheta\in \Tc_{\tau_{\zeta}}$ and let $(\vartheta_{n})_{n\ge 1}\subset  \Tc_{\vartheta}$ be such that $\vartheta_{n}\downarrow \vartheta$.  Since $\Yc^{\zeta_{}}_{\vartheta}=\Yc^{(\vartheta,X^{\zeta}_{\vartheta})}_{\vartheta}$, $\Yc^{\zeta_{}}_{\vartheta_{n}}=\Yc^{(\vartheta_{n},X^{\zeta}_{\vartheta_{n}})}_{\vartheta_{n}}$, and $X^{\zeta}_{\vartheta_{n}}=X^{(\vartheta,X^{\zeta}_{\vartheta})}_{\vartheta_{n}}$,   it follows from Corollary \ref{cor: stab time on the path of X} applied to the time intervalle $[\vartheta,\vartheta_{n}]$ that
\be\label{eq: proof right conti conv exp}
 \Yc^{\zeta}_{\vartheta}=\lim_{n\to \infty} \E_{\vartheta}[\Yc_{\vartheta_{n}}^{\zeta}].
 \ee
 We claim that $\limsup_{n\to \infty} \Yc_{\vartheta_{n}}^{\zeta}=\liminf_{n\to \infty} \Yc_{\vartheta_{n}}^{\zeta}$ a.s. Then, the above combined with  Remark \ref{rem: bound bar Y et Yc} and the dominated converge theorem  implies that
 $$
 \Yc^{\zeta}_{\vartheta}=\lim_{n\to \infty}  \Yc_{\vartheta_{n}}^{\zeta}.
 $$
 It remains to prove our claim. Let us first define
 $$
 \bar \eta:=\limsup_{n\to \infty} \Yc_{\vartheta_{n}}^{\zeta}\;\mbox{ and }\; \underline \eta:=\liminf_{n\to \infty} \Yc_{\vartheta_{n}}^{\zeta}.
 $$
 Assume on the contrary that the set
  \b*
 A:=\{\bar \eta>\underline \eta\} \; \mbox{ has positive probability.}
 \e*
 Note that  $A\in \Fc_{\vartheta}$, by right-continuity of the filtration.
 Then, define $\bar k_{1}=\underline k_{1}=1$ and
 \b*
 \bar k_{n+1}&:=&\min\{k> \bar k_{n}:  \Yc_{\vartheta_{k}}^{\zeta}\ge 2 \bar \eta/3 + \underline \eta/3\}\\
 \underline k_{n+1}&:=&\min\{k> \underline k_{n}:  \Yc_{\vartheta_{k}}^{\zeta}\le \bar \eta/3 +2\underline \eta/3\}
 \e*
 for $k\ge 1$, and set $\bar \vartheta_{n}:=\vartheta_{\bar k_{n}}\1_{A}+\vartheta_{{n}}\1_{A^{c}}$ and $\underline \vartheta_{n}:=\vartheta_{\underline k_{n}}\1_{A}+\vartheta_{{n}}\1_{A^{c}}$ for $n\ge 1$. It follows from the definition of $A$ that $\bar \vartheta_{n},\underline \vartheta_{n}$ are well-defined. They decrease to  $\vartheta$. Applying   \reff{eq: proof right conti conv exp} to the two sequences $(\bar \vartheta_{n})_{n\ge 1}$ and $(\underline \vartheta_{n})_{n\ge 1}$ and recalling the uniform bound of Remark \ref{rem: bound bar Y et Yc} leads to
\b*
(2 \bar \eta/3 + \underline \eta/3)\1_{A} \le \lim_{n\to \infty} \E_{\vartheta}[\Yc_{\bar \vartheta_{n}}^{\zeta}]\1_{A}
= \Yc^{\zeta}_{\vartheta}\1_{A}=\lim_{n\to \infty} \E_{\vartheta}[\Yc_{\underline \vartheta_{n}}^{\zeta}]\1_{A}\le ( \bar \eta/3 + 2\underline \eta/3)\1_{A},
\e*
a contradiction.

2. It follows from Proposition \ref{prop: dpp sup bar Y} that $\Yc^{\zeta}$ is a $f^{\zeta,0}$-supermartingale in the strong sense in the terminology of \cite{peng1999monotonic}. Since it is right-continuous, it follows from \cite[Theorem 3.3]{peng1999monotonic} that it admits left-limits. It is clear from Corollary \ref{cor: stab time on the path of X}    that it can not have  jumps on   $[0,T)$.

3. We finally prove the limit behavior at $T$.  It follows from Proposition \ref{prop: facelift}, Lemma \ref{lem: appendix linearisation}, \reff{eq: ass hat g} and the fact that $\Yc^{\zeta}$ is l\`ag  without positive jumps, see 1.~and 2.~above, that $ \Yc^{\zeta}_{T-}\ge \hat g(X^{\zeta}_{T})$.
  On the other hand, since $g(X^{\zeta}_{T})\le \hat g(X^{\zeta}_{T})$ and $ \hat g(X^{\zeta}_{T}+u)-\delta_{T}(u)\le \hat g(X^{\zeta}_{T})$ by construction and the fact that $\delta_{T}$ is sub-linear, we can  follow the arguments of the proof of Proposition \ref{prop: stab time} to deduce that
  \b*
  \Yc^{\zeta}_{\tau_{\zeta}\vee (T-\eps)}
  &\le&
  \E_{\tau_{\zeta}\vee (T-\eps)}\left[ \hat g(X^{\zeta}_{T})  \right]+C_{L}\;\eps^{\frac1{2}}.
\e*
By continuity of the filtration, this implies that
$
\Yc^{\zeta}_{T-}\le  \hat g(X^{\zeta}_{T}).
$
  \ep

\section{Weak versus strong formulation of the dual problem}\label{sec: weak versus strong}

The aim of this section is to prove that $\Yc^{\zeta}$ as defined in \reff{eq: def Yc} actually provides the minimal super-solution of \BSDE$_{K}(f,g,\zeta)$. Then, the statements of Theorem \ref{thm: regularity constrained bsde} will be a consequence of the results obtained in Section \ref{sec: estimates via the strong dual formulation}.  To this purpose, we first introduce the weak formulation associated to the optimal control problem  \reff{eq: def Yc}  and show that the value coincides. Then, we use standard arguments to show that this weak formulation actually provides the minimal super-solution of our constrained backward stochastic differential equation.

\subsection{Weak formulation}\label{sec: def weak form}

Given $\nu \in \Uc$ and $\zeta\in \Db_{2}$, we define the equivalent probability measure $\P^{\zeta,\nu}$ by
$$
\frac{d\P^{\zeta,\nu}}{d\P}= e^{-\frac12\int_{\tau_{\zeta}}^{T} |\sigma_{s}^{-1}(X^{\zeta}_{s})\nu_{s}|^{2} ds +\int_{\tau_{\zeta}}^{T}\sigma_{s}^{-1}(X^{\zeta}_{s})\nu_{s}dW_{s}}.
$$
Recall that $\sigma^{-1}$ is bounded by assumption.
Then,
\beq\label{eq: definition W nu}
W^{\zeta,\nu}:=W-\int_{\tau_{\zeta}}^{\tau_{\zeta}\vee \cdot}  \sigma_{s}^{-1}(X^{\zeta}_{s})\nu_{s} ds
\eeq
is a $\P^{\zeta,\nu}$-Brownian motion.

Given   $\vartheta\in \Tc_{\tau_{\zeta}}$   and  $G\in \Lb_{2}(\P^{\zeta,\nu},\Fc_{\vartheta})$, we set
$$
\Ecf{t,\vartheta}{\zeta,\nu}{G}:=U_{t\vee \vartheta}
$$
in which  $(U,V)\in \Sb_{2}(\P^{\zeta,\nu})\x \Hb_{2}(\P^{\zeta,\nu})$ satisfies
\be\label{eq: BSDE weak with G}
U_{t\vee \tau_{\zeta}}= G+\int_{t\vee \tau_{\zeta}}^{\vartheta} \left[f_{s}(X^{\zeta}_{s},U_{s},V_{s})-\delta_{s}(\nu_{s})\right] ds -\int_{t\vee \tau_{\zeta}}^{\vartheta} V_{s} dW^{\zeta,\nu}_{s},\;t\le T.
\ee
We finally define for $\tau \in \Tc_{\tau_{\zeta}}$
 \beq\label{eq: def Yb}
\Yb_{\tau}^{\zeta}:=\esssup\{\Ecf{\tau,T}{\zeta,\nu}{ g(X^{\zeta}_{T}) }:~\nu \in \Uc,\;\nu\1_{\[0,\tau\]}\equiv 0\} \;.
\eeq

\subsection{Equivalence of the strong and   weak formulations}

\begin{proposition}\label{prop: equiv strong and weak}
$
\Yb^{\zeta}_{\tau_{\zeta}}= \Yc^{\zeta}_{\tau_{\zeta}},
$
for  each $\zeta \in \Db_{2}$.
\end{proposition}

\proof We write $\zeta=(\tau,\xi)$. For sake of simplicity, we restrict to the situation $\tau=0$ so that $x_{0}:=\xi\in \R^{d}$. The general case is obtained by a conditioning argument. Let $\Uc^{\rm simple}$ denote the set of processes $\nu \in \Uc$ of the form
\be\label{eq: def Uc simple}
\nu=\sigma(X^{\zeta,\nu})\sum_{i=0}^{n-1} \phi_{i}\1_{(t_{i},t_{i+1}]}
\ee
in which $0=t_{0}<\cdots<t_{n}=T$ and $\phi_{i} \in \Lb_{\infty}(\Fc_{t_{i}})$ for $i\le n$. Note that $(\nu, X^{\zeta,\nu})$ is well-defined for any  $(\phi_{i})_{i\le n}\subset \Lb_{\infty}$ satisfying the previous measurability condition. This follows from \reff{eq: ass b sigma}.

 We define accordingly  $\tilde \Uc^{\rm simple}$ as the set of processes $\nu \in \Uc$ of the form
\be\label{eq: def tilde Uc simple}
\nu=\sigma(X^{\zeta})\sum_{i=0}^{n-1} \phi_{i}\1_{(t_{i},t_{i+1}]}
\ee
in which $0=t_{0}<\cdots<t_{n}=T$ and $\phi_{i} \in \Lb_{\infty}(\Fc_{t_{i}})$ for $i\le n$.

1. We first show that for each $\nu \in \Uc^{\rm simple}$ we can find $\tilde \nu  \in \Uc$ such that
$$
\Ec^{\zeta,\nu}_{\tau,T}[g(X^{\zeta,\nu}_{T})]= \Ecf{\tau,T}{\zeta,\tilde \nu}{ g(X^{\zeta}_{T}) }.
$$
Let $\nu \in \Uc^{\rm simple}$ be as in \reff{eq: def Uc simple} and note that we can identify $\phi_{i}$ to a Borel measurable map
$\omega \in \Omega\mapsto \phi_{i}(\omega)=\phi_{i}(\omega_{\cdot \wedge t_{i}})$, up to $\P$-null sets. Let us define $\tilde \nu$ by
$$
\tilde \nu(\omega)=\sigma(X^{\zeta})\phi_{i}\left(\omega^{\phi}_{\cdot \wedge t_{i}}\right) \;\mbox{ on } (t_{i},t_{i+1}]
$$
where $\omega^{\phi}$ is defined recursively by
$$
\omega_{s}^{\phi}:=\omega_{s}-\sum_{k=0}^{j-1} (t_{k+1}-t_{k})\phi_{k}(\omega^{\phi}_{\cdot \wedge t_{k}}) - (s-t_{j}) \phi_{j}(\omega^{\phi}_{\cdot \wedge t_{j}} ) \;\mbox{ for } s\in  (t_{j},t_{j+1}],
$$
with $\omega^{\phi}_{0}=0$.
Then, for $t\in (t_{i},t_{i+1}]$,
$$
X^{\zeta,\nu}_{t}=X^{\zeta,\nu}_{t_{i}}+\int_{t_{i}}^{t} (b_{s}(X^{\zeta,\nu}_{s})+\sigma_{s}(X^{\zeta,\nu}_{s})\phi_{i}(W_{\cdot\wedge t_{i}}))ds+\int_{t_{i}}^{t} \sigma_{s}(X^{\zeta,\nu}_{s})dW_{s}
$$
where $W$ is a Brownian motion under $\P$, while
$$
X^{\zeta}_{t}=X^{\zeta}_{t_{i}}+\int_{t_{i}}^{t} (b_{s}(X^{\zeta}_{s})+\sigma_{s}(X^{\zeta}_{s})\phi_{i}(W^{\tilde \nu}_{\cdot\wedge t_{i}}))ds+\int_{t_{i}}^{t} { \sigma_{s}(X^{\zeta}_{s})}dW^{\tilde \nu}_{s}
$$
where  $W^{\tilde \nu}$ is a Brownian motion under $\P^{\zeta,\tilde \nu}$.  This implies that the law of $(X^{\zeta,\nu},\nu)$ under $\P$ is the same as the law of $(X^{\zeta},\tilde \nu)$ under $\P^{\zeta,\tilde \nu}$. In view of Lemma \ref{lem : appendix schema num BSDE},   $\Ec^{\zeta,\nu}_{0,T}[g(X^{\zeta,\nu}_{T})]$ and $\Ecf{0,T}{\zeta,\tilde \nu}{ g(X^{\zeta}_{T}) }$ can be approximated by the same sequence of real numbers and are therefore equal.

2. The fact that   for each $\tilde \nu \in \tilde \Uc^{\rm simple}$ we can find $  \nu  \in \Uc$ such that
$$
\Ec^{\zeta,\nu}_{\tau,T}[g(X^{\zeta,\nu}_{T})]= \Ecf{\tau,T}{\zeta,\tilde \nu}{ g(X^{\zeta}_{T}) }
$$
follows from similar arguments.

3. To conclude the proof it remains to show that
\b*
\Yc_{\tau}^{\zeta}=\sup\{\Ec^{\zeta,\nu}_{\tau,T}[g(X^{\zeta,\nu}_{T})],\;\nu \in \Uc^{\rm simple}\}
\;\mbox{ and }\;
\Yb^{\zeta}_{\tau_{\zeta}}=\sup\{\Ecf{\tau,T}{\zeta,\tilde \nu}{ g(X^{\zeta}_{T}) },\;\tilde \nu \in \tilde \Uc^{\rm simple}\}.
\e*
We only prove the first identity, the second one being derived similarly. One inequality is trivial. Conversely, given any predictable and bounded process $\phi$,
we can find a bounded sequence of simple adapted processes $(\phi^{n})_{n\ge 1}$ such that $\E[\int_{0}^{T} |\phi^{n}_{s}-\phi_{s}|^{2}ds]\to 0$. By \reff{eq: ass b sigma}-\reff{eq: ass K borne},  $\nu^{n}:=\sigma(X^{\zeta,\nu^{n}})\phi^{n}\in \Uc$.  In particular, it follows from \reff{eq: ass b sigma} that $X^{\zeta,\nu_{n}}$ converges in $\Sb_{2}$ to $X^{\zeta,\nu}$ in which $\nu:=\sigma(X^{\zeta,\nu})\phi$.  Hence, after possibly passing to a subsequence,
$$
\liminf_{n\to \infty } g(X^{\zeta,\nu_{n}}_{T})\ge g(X^{\zeta,\nu}_{T})
$$
since $g$ is assumed to be lower-semicontinuous. By the comparison principle, Lemma \ref{lem: appendix stabilite et comp ?}, we have
$$
\Ec^{\zeta,\nu^{n}}_{\tau,T}[g(X^{\zeta,\nu^{n}}_{T})]\ge \Ec^{\zeta,\nu^{n}}_{\tau,T}[g(X^{\zeta,\nu}_{T})\wedge g(X^{\zeta,\nu_{n}}_{T})]
$$
in which $g(X^{\zeta,\nu}_{T})\wedge g(X^{\zeta,\nu_{n}}_{T})\to g(X^{\zeta,\nu}_{T})$ a.s. and in $\Lb_{2}$ by dominated convergence, recall \reff{eq: ass g}. We also have
$$
\E\left[\int_{\tau}^{T} |\delta_{s}(\nu_{s})-\delta_{s}(\nu^{n}_{s})|^{2}ds\right]\to 0
$$
since $(\delta_{t})_{t\le T}$ is equi-Lipschitz by \reff{eq: ass K borne}. Then, Lemma \ref{lem: appendix stabilite et comp ?} implies that
$$
\liminf_{n\to \infty }\Ec^{\zeta,\nu^{n}}_{\tau,T}[g(X^{\zeta,\nu^{n}}_{T})]\ge \lim_{n\to \infty} \Ec^{\zeta,\nu^{n}}_{\tau,T}[g(X^{\zeta,\nu}_{T})\wedge g(X^{\zeta,\nu_{n}}_{T})]  = \Ec^{\zeta,\nu}_{\tau,T}[g(X^{\zeta,\nu}_{T})].
$$
\ep
\subsection{Connection with the reflected backward stochastic differential equation}\label{sec: connection weak and K BSDE}

We now show that  $\Yb^{\zeta}$  identifies as the first component of the  minimal super-solution of the backward stochastic differential equation with constraint \BSDE$_{K}(f,g,\zeta)$.

\begin{Theorem}\label{thm: Doob Meyer decomposition}
  For all $\zeta \in \Db_{2}$,  there exists $\Zb^{\zeta}\in\Hb_{2}$ such that $(\Yb^{\zeta},\Zb^{\zeta})$ is the minimal supersolution of \BSDE$_{K}(f,g,\zeta)$.
\end{Theorem}

\proof The proof is standard  and written in the spirit of  \cite[Proof of Proposition 2.5]{CKS}. \\
1.   Similar arguments as in the proof of  Proposition \ref{prop: dpp sup bar Y} show that $\Yb$ satisfies a dynamic programming principle:   for all $\vartheta_{1}\le \vartheta_{2} \in \Tc$, such that $\tau_{\zeta}\le \vartheta_{1}$, we have
\be\label{eq: dpp for Ycweak}
\Yb_{\vartheta_{1}}^{\zeta}=\esssup_{\nu \in \Uc} \Ecf{\vartheta_{1},\vartheta_{2}}{\zeta,\nu}{\Yb_{\vartheta_{2}}^{\zeta} }.
\ee
We also observe that $\Yb^{\zeta}$ is c\`ad. This follows from Proposition \ref{prop: equiv strong and weak} and  Proposition \ref{prop: path continuity}. Then, the non linear Doob-Meyer decomposition of \cite[Theorem 3.3]{peng1999monotonic} implies the existence of $\Zb^{\zeta,\nu}\in \Hb_{2}(\P^{\nu,\zeta})$ and of a c\`adl\`ag non-decreasing adapted process $\Kb^{\zeta,\nu}$ such that
$$
\Yb^{\zeta}_{\vartheta} =   g(X^{\zeta}_{T})+\int_{\vartheta}^{T} \left(f_{s}(X^{\zeta}_{s},\Yb^{\zeta}_{s},\Zb^{\zeta,\nu}_{s})-\delta_s(\nu_{s})\right)ds-\int_{\vartheta}^{T}\Zb^{\zeta,\nu}_{s}dW^{\zeta,\nu}_{s}+\Kb^{\zeta,\nu}_{T}-\Kb^{\zeta,\nu}_{\vartheta} \;,  \; \vartheta\ge \tau_{\zeta}\;.
$$
By identification of the It\^{o} decomposition under each $\P^{\zeta,\nu}$, we obtain $\Zb^{\zeta,\nu}=\Zb^{\zeta,0}=:\Zb^{\zeta}$. Moreover,  \reff{eq: definition W nu} implies that for any $\nu\in\Uc$ and $\vartheta\in\Tc_{\tau_{\zeta}}$
$$
\Kb^{\zeta,0}_{T}-\Kb^{\zeta,0}_{\vartheta}=\Kb^{\zeta,\nu}_{T}-\Kb^{\zeta,\nu}_{\vartheta}+\int_{\vartheta}^{T} (\Zb^{\zeta}_{s}\sigma^{-1}(X^{\zeta}_{s}) \nu_{s}-\delta_s(\nu_{s})) ds
\ge \int_{\vartheta}^{T} (\Zb^{\zeta}_{s}\sigma^{-1}(X^{\zeta}_{s}) \nu_{s}-\delta_s(\nu_{s})) ds\;,
$$
 since $\Kb^{\zeta,\nu}$ is non-decreasing. By using a similar measurable selection argument as in the proof of Proposition \ref{prop: facelift}, this shows that
 $$
 \inf_{|u|=1}(\delta(u)-\Zb^{\zeta}\sigma^{-1}(X^{\zeta}) u)\ge 0\;\; dt\otimes d\P\mbox{-a.e.}
\;\Leftrightarrow\;
 \Zb^{\zeta}\sigma^{-1}(X^{\zeta}) \in K\;\;dt\otimes d\P\mbox{-a.e.},
 $$
 see e.g. \cite{Rockafellar}.
Since $\Kb^{\zeta,0}$ is non-decreasing, writing the above for $\nu=0$ implies that $(\Yb^{\zeta},\Zb^{\zeta})$ is a super-solution of \BSDE$_{K}(f,g,\zeta)$.

2. We now prove the minimality property. If $(U,V)$ is a super-solution of  \BSDE$_{K}(f,g,\zeta)$, then it follows from the definition of $\delta$ and \reff{eq: definition W nu} that it is also a super-solution of \reff{eq: BSDE weak with G} with $G=g(X^{\zeta}_{T})$, for each $\nu\in \Uc$. In particular, $U\ge  \Ecf{\cdot,T}{\zeta,\nu}{ g(X^{\zeta}_{T}) }$, and we conclude by arbitrariness of $\nu \in \Uc$.
\ep

\section{Possible extensions}\label{sec: possible extensions}

In order to focus on the main ideas, we have restricted ourselves to a rather stringent framework. Some of the conditions used in this paper can certainly be weakened on a case by case basis. We discuss here some  straightforward extensions or variations.

\subsection{Invertibility condition}

We have assumed that $\sigma$ is invertible but all our arguments go through if we add a component $ X^{o}$ to $X$ which has a dynamic of the form
$$
d  X^{o}_{t}= b^{o}_{t}( X^{o}_{t},X_{t})dt
$$
with $ b^{o}$ Lipschitz and bounded in space, uniformly in time. Then $X^{\zeta,\nu}$ has to be replaced by $\bar X^{\zeta,\nu}=( X^{o,\zeta,\nu},X^{\zeta,\nu})$ with dynamics
\b*
d\bar X^{\zeta,\nu}_{t}=\left(\begin{array}{c} b^{o}_{t}(\bar X^{\zeta,\nu}_{t}) \\  b_{t}( \bar X^{\zeta,\nu}_{t})+\nu_{t} \end{array}\right) dt+
\left(\begin{array}{c} 0\\  \sigma_{t}(\bar X^{\zeta,\nu}_{t}) \end{array}\right) dW_{t}.
\e*
The face-lift of $g$ is defined accordingly
$$
\hat g(x^{o},x):=\sup_{u\in \R^{d}} \left(g(x^{o},x+u)-\delta_{T}(u)\right),
$$
and so on.

The case of a general  non-invertible coefficient $\sigma$ can be treated along the lines of \cite{bouchard2010obstacle}, in which it is explained how the face-lift should then be  performed.

\subsection{Direct constraint on the gains-process}

The contraint \reff{eq: def constrained bsde 2} is motivated by financial applications in which the component $V$ can be interpreted as the number of risky assets $X$ held in an hedging portfolio for the contingent claim $g(X^{\zeta}_{T})$, see \cite{Broadie}. It can be replaced by
\beq\label{eq: def constrained bsde on Z}
V \in K\;\;dt\otimes d\P{\rm -a.e} \; \mbox{ on } [\![\tau_{\zeta},T]\!],
\eeq
when $\sigma$ does not depend on $x$.

In this case, $X^{\zeta,\nu}$  must be taken of the form
\b*
d  X^{\zeta,\nu}_{t}= \left( b_{t}(  X^{\zeta,\nu}_{t})+\sigma_{t}\nu_{t} \right) dt+\sigma_{t} dW_{t}.
\e*
If one assumes that $t\mapsto \sigma_{t}$ is right-continuous on $[0,T)$ and left-continuous at $T$, then \reff{eq: ass hat g}-\reff{eq: face lift t<T} become
\b*
\hat g(x)&=&\sup_{u\in \R^{d}} \left(g(x+\sigma_{T}u)-\delta_{T}(u)\right)
\\
 \Yc^{\zeta}_{\tau_{\zeta}}&=& \esssup_{u\in \Lb_{\infty}(\R^{d},\Fc_{\tau_{\zeta}})}(\Yc^{(\tau_{\zeta},{\xi_{\zeta}+\sigma_{\tau_{\zeta}}u})}_{\tau_{\zeta}}-\delta_{{\tau_{\zeta}}}(u))\;\mbox{ a.s. on $\{\tau_{\zeta}<T\}$.}
\e*
The change of measure for the weak formulation is
$$
\frac{d\P^{\zeta,\nu}}{d\P}= e^{-\frac12\int_{\tau_{\zeta}}^{T} |\nu_{s}|^{2} ds +\int_{\tau_{\zeta}}^{T}\nu_{s}dW_{s}}.
$$
In particular, we do not need $\sigma$ to be invertible anymore.

The results of Theorem \ref{thm: regularity constrained bsde} are obtained by following step by step the arguments used in this paper up to the modifications described above.

Moreover, the boundedness condition \reff{eq: ass K borne} can then be weakened. Indeed, it can be avoided by using \reff{eq : ass delta continuous}-\reff{eq : ass delta non-increasing} in all our proofs, except in the proof of Proposition \ref{prop: equiv strong and weak} in which it is used twice. First to ensure that the controls $\nu$ constructed from the families $(\phi_{i})_{i\le n}$ satisfy $\delta(\nu)<\infty$.  But in the case \reff{eq: def constrained bsde on Z}, the $\nu$'s are of the form $\sum_{i=0}^{n-1}\phi_{i}\1_{(t_{i},t_{i+1}]}$. Taking $\delta_{t_{i}}(\phi_{i})<\infty$ is then enough. It is also used in the approximation argument of Step 3, as it implies that  $(\delta_{t})_{t\le T}$ is equi-Lipschitz, but for the constraint \reff{eq: def constrained bsde on Z},  it suffices to assume, for instance, that the domain of $\delta_{t}$ does not depend on $t$, which means that the directions in which $K_{t}$ is bounded do not depend on $t$.  This is not enough when $\nu$ is of the form used in the proof of Proposition \ref{prop: equiv strong and weak} because of the transformation through the  matrix $\sigma$, unless additional assumptions are made on it.

 Our arguments are not valid if $\sigma$ depends on $x$ because the coefficients driving  $X^{\zeta,\nu}$ are no more  Lipschitz uniformly in the control. This is crucial for Corollary \ref{cor: space regularity}.

\subsection{Optimal control of constrained BSDEs}

One can allow the coefficients $b,\sigma$ and $f$ to depend on an additional control $\alpha$ in a set $\Ac$ of predictable processes with values in a compact set $A\subset \R^{d}$. Then, all our proofs go through whenever the conditions \reff{eq: ass b sigma}-\reff{eq: ass f} are uniform with respect to this additional control, and the coefficients are continuous in this additional variable. The arguments used in Section \ref{sec: estimates via the strong dual formulation} do not change. It is the same for Section \ref{sec: connection weak and K BSDE}, for $\alpha\in \Ac$ given.  However, a  continuity assumption on the coefficients with respect to the control will be required to prove the counterpart of Proposition \ref{prop: equiv strong and weak}: the approximation by step constant processes has to be applied  to $(\nu,\alpha)$ in place of $\nu$.

\subsection{Random coefficients with delay}

One can also assume that the coefficients $b,\sigma$ and $f$ are random, satisfying the usual predictability condition,  whenever the conditions \reff{eq: ass b sigma}-\reff{eq: ass f} are uniform in $\omega$. Again, the arguments of Section \ref{sec: estimates via the strong dual formulation} and  Section \ref{sec: connection weak and K BSDE} do not change. However, the proof of Proposition \ref{prop: equiv strong and weak} can not be adapted unless the dependence holds with a fixed delay: there exists $\iota>0$ such that, for all $t\le T$,  $b_{t},\sigma_{t},f_{t}$ depends on $\omega$ only through $(\omega_{s})_{s\le t-\iota}$.  With this condition, Steps 1.~and 2.~remain correct for simple processes associated to a time grid $\{t_{i},i\le n\}$  such that $\max\{t_{i+1}-t_{i},i\le n-1\}\le \iota$.  {As in the optimal control case, Step 3 also requires some continuity of the coefficient in $\omega$, e.g.~uniform continuity for the usual sup-norm topology. }

\appendix
\section{Auxiliary results}

We collect here some standard results that have been used all over this paper.

In this section, we denote by $\Dc_{b}$ the set of  measurable maps $\psi: \Omega\x [0,T]\x \R^{d}\x \R\mapsto \R$  such that
$(\psi_{t}(y,z))_{t\le T}$ is progressively measurable for all $(y,z) \in \R\x \R^{d}\mapsto \R$ and
$$
|\psi(y,z)-\psi(0)|\le L_{\psi}(|y|+|z|) \mbox{ for all } (y,z)\in \R\x \R^{d},~dt\otimes d\P-a.e.
$$
for some constant $L_{\psi}>0$. Given $(\vartheta, G)\in \Db_{2}$ and $\tau\in \Tc$ such that $\tau \le \vartheta$, we set $\Ec^{\psi}_{\tau,\vartheta}(G):=U_{\tau}$ where $(U,V)\in \Sb_{2}\x \Hb_{2}$ is the solution of
\be\label{eq: def BSDE appendix}
U_{t\vee \tau}= G+\int_{t\vee \tau}^{\vartheta} \psi_{s}(U_{s},V_{s})ds -\int_{t\vee \tau}^{\vartheta} V_{s} dW_{s}, \;t\in [0,T].
\ee

\begin{Lemma}\label{lem: appendix linearisation}  Fix $\psi \in \Dc_{b}$ and  $(\vartheta, G)\in \Db_{2}$. Then, for all $\tau \in \Tc$ such that $\tau \le \vartheta$:
\begin{enumerate}[{\rm (a)}]
\item We have
\b*
\Ec^{\psi}_{\tau,\vartheta}(G)= \E_{\tau}[H^{\tau }_{\vartheta}G+\int_{\tau}^{\vartheta} H^{\tau }_{s} \psi_{s}(0) ds ]
\e*
where $H^{\tau }$ solves
$$
H^{\tau }=1+\int_{\tau}^{\cdot}  \kappa^{Y}_{t} H^{\tau }_{t}dt + \int_{\tau}^{\cdot}   \kappa^{Z}_{t} H^{\tau }_{t}dW_{t}
$$
for some predictable processes $\kappa^{Y}$ and $\kappa^{Z}$ that are bounded by a constant which only depends on $L_{\psi}$.

In particular, if there exists a constant $c>0$ such that $\E[|G|^{2}]\le c$ and $|\psi(0)|\le c$ $dt\otimes d\P$, then
\b*
|\Ec^{\psi}_{\tau,\vartheta}(G)-\E_{\tau}[G]|\le C\,\E_{\tau}[\vartheta-\tau]^{\frac12},
\e*
for some $C>0$ which depends only on $c$ and $L_{\psi}$. 

 \item If $(U,V)\in \Sb_{2}\x \Hb_{2}$ satisfies
$$
U_{t\vee \tau}\le  G+\int_{t\vee \tau}^{\vartheta} \psi_{s}(U_{s},V_{s})ds -\int_{t\vee \tau}^{\vartheta} V_{s} dW_{s}, \;t\in [0,T],
$$
then
\b*
U_{\tau}\le  \E_{\tau}[H^{\tau }_{\vartheta}G+\int_{\tau}^{\vartheta} H^{\tau }_{s} \psi_{s}(0) ds ]
\e*
where $H^{\tau}$ is defined as in {\rm (a)}.
\end{enumerate} 
\end{Lemma}

\proof This follows from a standard linearization argument, see e.g.~\cite[Proof of Theorem 1.6]{Par98}.
\ep

\begin{Lemma}\label{lem: appendix stabilite et comp ?}  Fix $\psi^{1},\psi^{2} \in \Dc_{b}$ and  $(\vartheta, G_{1}),(\vartheta, G_{2})\in \Db_{2}$.
\begin{enumerate}[{\rm (a)}]
\item Assume that there exists a process $\kappa$ such that
$$
  |\psi^{1}-\psi^{2}|(\thetar)\le \kappa\;\; \;dt\otimes d\P
$$
for all $\thetar\in \R\times \R^{d}$.
 Then, for all $\tau \in \Tc$ such that $\tau \le \vartheta$,
\b*
|\Ec^{\psi^{1}}_{\tau,\vartheta}(G_{1})- \Ec^{\psi^{2}}_{\tau,\vartheta}(G_{2})|\le C \E_{\tau}\left[|G_{1}-G_{2}|^{2}+\int_{\tau}^{T}|\kappa_{s}|^{2} ds  \right]^{\frac12}
\e*
where $C>0$ is a constant which  depends only on $L_{\psi^{1}}$ and $L_{\psi^{2}}$.

\item Assume that $G_{1}\le G_{2}$ and
$
\psi^{1}(\thetar)\le \psi^{2}(\thetar)$ $dt\otimes d\P
$
for all $\thetar\in \R\times \R^{d}$. Then, $\Ec^{\psi^{1}}_{\tau,\vartheta}(G^{1})\le \Ec^{\psi^{2}}_{\tau,\vartheta}(G^{2}) $ for all $\tau \in \Tc$ such that $\tau \le \vartheta$.
\end{enumerate} 
\end{Lemma}

\proof The first assertion follows from
\cite[Theorem 1.5]{Par98}. The second one is \cite[Theorem 1.6]{Par98}. \ep

\begin{Lemma}\label{lem: appendix liminf eps and E condieps} Let $(G_{\eps})_{\eps>0}$ be a family of random variable, uniformly bounded in $\Lb_{1}$, and let $G\in \Lb_{1}$
be such that $\liminf_{\eps\to 0} G_{\eps}\ge G$. Let $(\tau_{\eps})_{\eps>0}$ be a sequence of stopping times such that $\lim_{\eps\to 0} \tau_{\eps}=\tau \in \Tc$. Then, there exists a  sequence  $(\eps_{n})_{n\ge 1}\subset (0,1)$ such  that
$$
\liminf_{n\to \infty} \E_{\tau_{\eps_{n}}}[G_{\eps_{n}}]\ge \E_{\tau}[G] \mbox{ and } \lim_{n\to \infty} \eps_{n}=0.
$$
\end{Lemma}

\proof We write
$$
 \E_{\tau_{\eps}}[G_{\eps}]=\E_{\tau_{\eps}}[G] + \E_{\tau_{\eps}}[G_{\eps}-G]\ge  \E_{\tau_{\eps}}[G] - \E_{\tau_{\eps}}[(G_{\eps}-G)^{-}].
$$
The first term on the right-hand side converges a.s.~to   $\E_\tau[G]$ by the continuity of the martingales in a Brownian filtration.  The second term converges in $\Lb_{1}$ to $0$, and therefore a.s.~along a subsequence.
\ep

\begin{Lemma}\label{lem : appendix schema num BSDE} Fix $\psi\in  \Dc_{b},G\in \Lb_{2}(\Fc_{T})$. Set $t_{i}^{n}:=iT/n$ for $i\le n$, $n\ge 1$, and define recursively for  $i=n-1,\ldots,0$
 \b*
 U^n_{t^{n}_{i}}&=&\E_{t^{n}_{i}}[ U^n_{t^{n}_{i+1}}+\int_{t^{n}_{i}}^{t^{n}_{i+1}} \psi_{s}( U^n_{t^{n}_{i}}, V^n_{t^{n}_{i}})ds]\;,\;
 V^n_{t^{n}_{i}}=(t^{n}_{i+1}-t^{n}_{i})^{-1}\E_{t^{n}_{i}}[U^n_{t^{n}_{i+1}}(W_{t^{n}_{i+1}}-W_{t^{n}_{i}})]
 \e*
 in which $U^{n}_{T}:=G$. Then, $U^{n}_{0}\to \Ec^{\psi}_{0,T}[G]$ as $n\to \infty$.
\end{Lemma}

\proof It suffices to repeat the argument of \cite[Proof of Theorem 3.1]{bouchard2004discrete} and observe that their estimate contained in  \cite[Lemma 3.2]{bouchard2004discrete} is not needed if we are not interested by the speed of convergence. Indeed, one can simply  use the fact that, if $(U,V)$ denotes the solution of  \reff{eq: def BSDE appendix} with $\tau=0$, then
$$
\max_{1\le i\le n} \E[\sup_{t^{n}_{i-1}\le t\le t^{n}_{i}} |U_{t}-U_{t^{n}_{i-1}}|^{2}] + \sum_{i=1}^{n} \E[\int_{t^{n}_{i-1}}^{t^{n}_{i}}|V_{t}-\bar V^{n}_{t^{n}_{i-1}}|^{2}dt]\to 0,
$$
in which
$$
\bar V^{n}_{t^{n}_{i-1}}:=(t^{n}_{i}-t^{n}_{i-1})^{-1}\E_{t^{n}_{i-1}}[\int_{t^{n}_{i-1}}^{t^{n}_{i}} V_{t} dt ].
$$
The convergence of the left-hand side term is standard, it follows from the continuity of the path of $U$ which belongs to $\Sb_{2}$.  The convergence of the second term is also clear since $\bar V^{n}$ provides the best approximation of $V$ in $\Lb_{2}(dt\otimes d\P)$ by a step constant process on $\{t^{n}_{i},i\le n\}$.
\ep

\bibliographystyle{plain}

\end{document}